\documentclass[11pt,leqno]{article}
\usepackage{amsfonts}
\usepackage{amsmath}
\usepackage{amssymb}
\newcommand{\pref}[1]{(\ref{#1})}
\newtheorem{theo}{Theorem}[section]
\newtheorem{lem}[theo]{Lemma}
\newtheorem{prop}[theo]{Proposition}

\title{\Large \bf {Hermitian spin surfaces with small eigenvalues of the Dolbeault operator}}
\author{{\sc Bogdan Alexandrov}
\thanks{Supported by SFB 288 "Differential geometry and quantum
physics" and SPP 1154 "Globale Differentialgeometrie" of DFG	and The European Contract Human Potential Programme,
Research Training Network HPRN-CT-2000-00101}
}
\date{}
\begin{document}
\maketitle
\vspace{5mm}
\begin{abstract}
We study the compact Hermitian spin surfaces with positive conformal scalar curvature on which the first eigenvalue of the Dolbeault operator of the spin structure is the smallest possible. We prove that such a surface is either a ruled surface or a Hopf surface. We give a complete classification of the ruled surfaces with this property. For the Hopf surfaces we obtain a partial classification and some examples.
\\[10mm]
\\[10mm]
\end{abstract}

\section{Introduction}

It was proved by Friedrich in \cite{F1} that on a compact $n$-dimensional Riemannian spin manifold $M$ with positive scalar curvature $s$ any eigenvalue $\lambda$ of the Dirac operator satisfies the inequality
\begin{equation}\label{1}
\lambda ^2 \ge \frac{n}{n-1} \frac{{\rm min}_M \, s}{4}.
\end{equation}
The limiting manifolds, i.e., the manifolds on which the first eigenvalue satisfies the equality in \pref{1}, are Einstein and have holonomy $SO(n)$. In the simply connected case their classification has been completed by B\"ar \cite{Ba}.

The fact that the holonomy of a limiting manifold is $SO(n)$ implies, for example, that \pref{1} is strict on K\"ahler manifolds. For K\"ahler manifolds a better estimate was found by Kirchberg \cite{K1,K2}. In the 4-dimensional case (this is the dimension we are interested in in this paper) it reads as follows:
\begin{equation}\label{2}
\lambda ^2 \ge \frac{{\rm min}_M \, s}{2}.
\end{equation}
The limiting manifolds are characterized by having constant scalar curvature and the existence of a K\"ahlerian twistor spinor \cite{K2}. It was shown by Friedrich \cite{F2} that they are biholomorphically equivalent to $\mathbb{C}P^1 \times \mathbb{C}P^1$ or $T
\times \mathbb{C}P^1$, where $T$ is a torus.

Now a natural question arises: Is it possible to generalize in a reasonable way these results on Hermitian (non-K\"ahler) manifolds?

In \cite{Hit} Hitchin showed that there is a bijective correspondence between the spin structures on a Hermitian manifold and the holomorphic square roots $S$ of its canonical bundle $K$. If the manifold is moreover K\"ahler, then the Dirac operator coincides with the Dolbeault operator $\Box =\sqrt 2 (\overline {\partial}_S + {\overline {\partial}}^*_S)$ of the corresponding square root $S$.

This shows that in order to capture better the existence of a complex structure it would be perhaps reasonable to replace the Dirac operator on a Hermitian manifold by the Dolbeault operator of $S$.

Indeed, the author together with Grantcharov and Ivanov proved in \cite{AGI} the following

\begin{theo}\label{th1}
Let M be a compact Hermitian spin surface with positive conformal scalar curvature $k$. Then for each eigenvalue $\lambda$ of $\Box$
\begin{equation}\label{3}
\lambda ^2 \ge \frac{{\rm min}_M \, k}{2} .
\end{equation}
The first eigenvalue satisfies the equality in \pref{3} iff $k$ is constant and there exists a non-identically zero section  $\psi \in \Gamma(S^*)$ such that
\begin{equation}\label{305}
(\nabla^{-3})^{0,1}\psi = 0.
\end{equation}
In this case $M$ is locally conformally K\"ahler.
\end{theo}

The connection $\nabla^{-3}$ is a member of the one-parameter family of canonical Hermitian connections $\nabla ^t$, $t \in {\bf R}$ (see \cite{G2}), and $(\nabla^{-3})^{0,1}$ is the $(0,1)$-part of $\nabla^{-3}$. On a K\"ahler manifold $\nabla ^t$ coincides with the Levi-Civita connection, the conformal scalar curvature coincides with the scalar curvature and a solution of \pref{305} is just a holomorphic section of $S^*$, i.e., a K\"ahlerian twistor spinor. Thus Theorem~\ref{th1} generalizes the above cited result of Kirchberg.

In this paper we investigate the Hermitian surfaces on which the limiting case in Theorem~\ref{th1} occurs. From now on we call them {\it limiting surfaces}.

In section 2 we prove that such a manifold is either a ruled surface or a Hopf surface. In the next two sections we discuss these two cases respectively.

A ruled surface of genus $g$ is a holomorphic $\mathbb{C}P^1$-bundle over a complex curve of genus $g$. It admits a K\"ahler metric and therefore every locally conformally K\"ahler metric is globally conformally K\"ahler by a result of Vaisman \cite{V1}. Thus, since the existence of a non-zero solution of \pref{305} is a conformally invariant condition, it is equivalent to the existence of a square root $S$ of the canonical line bundle such that $S^*$ has holomorphic section. This property characterizes the limiting ruled surfaces as the existence of a Hermitian metric with positive constant conformal scalar curvature is ensured by the existence of a K\"ahler metric with positive scalar curvature on each ruled surface, the latter being a special case of a theorem of Yau \cite{Y}. The limiting ruled surfaces, together with the corresponding spin structures, are described in Theorem~\ref{th2}. There are limiting ruled surfaces of arbitrary genus. In particular, there are many more of them than in the K\"ahler case. Even for genus $g=0$ or $1$ there are limiting surfaces other than the trivial bundles $\mathbb{C}P^1 \times \mathbb{C}P^1$ and $T\times \mathbb{C}P^1$. For example, all even Hirzebruch surfaces are limiting. It is also interesting to remark that some ruled surfaces are limiting with only one of their spin structures, while others are limiting with all of them.

A Hopf surface is a complex surface whose universal cover is $\mathbb{C}^2 \backslash 0$. Every Hopf surface is finitely covered by $S^1 \times S^3$. Thus the first Betti number is 1  and therefore there exists no K\"ahler metric. Nevertheless we can describe the Hopf surfaces on which the 'holomorphic' condition \pref{305} is satisfied. The description involves the cohomology class of the Lee form of the metric and is given in Theorem~\ref{th3}. This theorem would give a classification of the limiting Hopf surfaces if the possible Lee forms of locally conformally K\"ahler metrics with positive conformal scalar curvature were known. There are some restrictions on Hopf surfaces of class 1 coming from the fact that they admit Vaisman metrics \cite{GO,Bel} and a result of Tsukada \cite{Ts}. But in general it is not even known whether each Hopf surface admits a locally conformally K\"ahler metric with positive conformal scalar curvature. Still Theorem~\ref{th3} allows us to obtain numerous examples of limiting Hopf surfaces and also of such Hopf surfaces which cannot be limiting.

\section{Preliminaries}

Let $(M,h,I)$ be a $2m$-dimensional Hermitian manifold, i.e.,  $I$ is an (integrable) complex structure  and $h$ is a Riemannian metric such that $h(I\cdot,I\cdot)=h$. We denote by $\Omega$ the K\"ahler form,
$$\Omega(X,Y) = h(X,IY).$$
Let $\Lambda^{p,q} M = \Lambda^p (T^{1,0})^* \otimes \Lambda^q (T^{0,1})^*$ be the bundle of forms of type $(p,q)$. The canonical line bundle is $K=\Lambda ^{m,0} M$.

It was proved in \cite{Hit} that the spin structures on a Hermitian manifold are in one-to-one correspondence with the set of holomorphic line bundles $S$ such that $S \otimes S = K$. Given a spin structure $S$, the corresponding spinor bundle is
\begin{equation}\label{4}
\Sigma M =  \sum _{r=0}^m \Lambda ^{0,r} M \otimes S.
\end{equation}
We denote by $\Box$ the Dolbeault operator of $S$,
$$\Box =\sqrt 2 (\overline {\partial}_S + {\overline {\partial}}^*_S).$$
Thus the Dirac operator $D$ and the Dolbeault operator $\Box$ act on sections of $\Sigma M$ (in particular, on sections of $\Lambda ^{0,m} M \otimes S = S^*$) and if the manifold is K\"ahler, they coincide.

Recall that a Hermitian connection on a Hermitian manifold is a connection with respect to which both the metric $h$ and the complex structure $I$ are parallel. In particular, such a connection defines a connection in the spinor bundle $\Sigma M$ which preserves the decomposition \pref{4}. There is a 1-parameter family $\nabla ^t $, $t \in {\bf R}$, of distinguished Hermitian connections, the canonical Hermitian connections (see \cite{G2}). They are defined as follows:
$$\nabla ^t = (1-t)\nabla ^0 + t\nabla^1, \qquad t \in {\bf R},$$
where $\nabla ^0$ is the projection of the Levi-Civita connection $\nabla$ in the affine space of all Hermitian connections and $\nabla ^1$ is the Chern connection, i.e., the unique Hermitian connection such that when considered as a connection on the holomorphic tangent bundle $T^{1,0}M$ its $(0,1)$-part $(\nabla^1)^{0,1} = \overline {\partial}_{T^{1,0}M}$. On K\"ahler manifolds the canonical Hermitian connections coincide with the Levi-Civita connection.

Although the following notions can be defined in arbitrary dimensions, from now on we restrict our considerations to Hermitian surfaces, i.e., to complex dimension $m=2$.

In this case the wedge product with the K\"ahler form $\Omega$ yields an isomorphism between  the spaces of 1-forms and 3-forms. The Lee form of $(M,h,I)$ is defined to be the unique 1-form $\theta$ such that
$$d\Omega = \theta \wedge \Omega.$$

The conformal scalar curvature $k$ is the scalar curvature with respect to $h$ of the canonical Weyl connection of $(M,h,I)$. For our purposes it will be enough to define it through the explicit formula
\begin{equation}\label{405}
k = s - \frac {3}{2}|\theta |^2 -3d^* \theta ,
\end{equation}
where $s$ is the scalar curvature of $h$.

If we change the metric conformally, $\widetilde{h} = e^f h$, then
$$\widetilde{\theta }= \theta + df, \qquad \widetilde{k}=e^{-f} k.$$
So we get the following straightforward consequences:
\begin{itemize}
\item $(M,h,I)$ is K\"ahler iff $\theta =0$ and in this case $k=s$.
\item $(M,h,I)$ is globally conformally K\"ahler iff $\theta $ is exact.
\item $(M,h,I)$ is locally conformally K\"ahler iff $\theta $  is closed.
\item If $k>0$, then there exists a conformally equivalent metric with positive constant conformal scalar curvature: if $c>0$, then $\widetilde{h} = \frac{k}{c} h$ has $\widetilde{k}= c$.
\end{itemize}

The next straightforward lemma (cf Lemma~4 in \cite{AGI}) shows that the existence of a non-zero solution of \pref{305} is a conformally invariant property.

\begin{lem}\label{lem1}
If $\psi \in \Gamma(S^*)$ is a solution of \pref{305} for the metric $h$, then $e^f \psi$ is a solution of \pref{305} for $\widetilde{h} = e^f h$.
\end{lem}

Thus, to find the limiting surfaces it will be enough to find the Hermitian surfaces with positive (not necessarily constant) conformal scalar curvature admitting a non-zero solution of \pref{305}.

Let us consider the exact sequence
$$0 \longrightarrow H^1 (M,\mathbb{Z}) \longrightarrow H^1 (M,\mathbb{C}) \longrightarrow H^1 (M,\mathbb{C}^*) \longrightarrow \dots$$
coming from the exact sequence
$$0 \longrightarrow \mathbb{Z} \longrightarrow \mathbb{C} \stackrel{e^{2\pi i \bullet}}{\longrightarrow}  \mathbb{C}^* \longrightarrow 0.$$
For $[\omega] \in H^1 (M,\mathbb{C})$ we denote by $E([\omega])$ its image in $H^1 (M,\mathbb{C}^*)$. Thus $E([\omega])$ is a complex line bundle with constant transition functions and therefore holomorphic. Let $\{U_\alpha\}$ be a good open cover of $M$ and $\omega |_{U_\alpha} = df_\alpha$. Then $f_\alpha - f_\beta$ are constant on $U_\alpha \cap U_\beta$ and the transition functions of $E([\omega])$ with respect to $\{U_\alpha\}$ are $e^{2\pi i (f_\alpha - f_\beta)}$.

\begin{prop}\label{prop0}
There exists a non-zero solution of \pref{305} on a compact locally conformally K\"ahler spin surface $(M,h,I,S)$ with Lee form $\theta$ iff $H^0 (M,\mathcal{O}(E([ - \frac{1}{2\pi i}\theta]) \otimes S^*)) \not = 0$.
\end{prop}

\noindent {\it Proof:}
Let $\{U_\alpha\}$ be a good open cover of $M$, $\theta |_{U_\alpha} = df_\alpha$, $s_\alpha$ be non-zero holomorphic sections of $S^*|_{U_\alpha}$. Let $\psi$ be a solution of \pref{305}, $\psi|_{U_\alpha} = \psi_\alpha s_\alpha$. Then $\psi_\alpha = g_{\alpha \beta} \psi_\beta$, where $g_{\alpha \beta}$ are the transition functions of $S^*$ determined by $\{s_\alpha \}$.

Since $e^{-f_\alpha} h$ is a K\"ahler metric on $U_\alpha$, Lemma~\ref{lem1} yields that $e^{-f_\alpha} \psi|_{U_\alpha} = e^{-f_\alpha} \psi_\alpha s_\alpha$ is a holomorphic section of $S^*$ over $U_\alpha$. Hence the functions $\varphi_\alpha := e^{-f_\alpha} \psi_\alpha$ are holomorphic. The transition functions of $E([-\frac{1}{2\pi i}\theta]) \otimes S^*$ are $h_{\alpha \beta} = e^{-f_\alpha + f_\beta} g_{\alpha \beta}$ and obviously $\varphi_\alpha = h_{\alpha \beta} \varphi_\beta$. Thus $\{\varphi_\alpha \}$ define a holomorphic section of $E([-\frac{1}{2\pi i}\theta]) \otimes S^*$, i.e., $H^0 (M,\mathcal{O}(E([-\frac{1}{2\pi i}\theta]) \otimes S^*)) \not = 0$.

The converse is proved in the converse way. \hfill $\Box $

\begin{prop}\label{prop1}
A compact spin Hermitian surface $(M,h,I,S)$ with positive conformal scalar curvature is biholomorphically equivalent to a ruled surface or a Hopf surface.
\end{prop}

\noindent {\it Proof:}
Since $M$ is spin, $(M,I)$ is a minimal complex surface. Indeed, if $C$ is an exceptional curve, then $(c_1 (C) \wedge c_1 (K))[M] = -1$ which is impossible because $c_1 (K) = 2 c_1 (S)$.
$M$ is compact and therefore there exists a Gauduchon metric in the conformal class of $h$, i.e., a metric whose Lee form is co-closed (see \cite{G3}). The positivity of the conformal scalar curvature is conformally invariant, so we can assume that $h$ is the Gauduchon metric. Now \pref{405} shows that the scalar curvature of $h$ is also positive. Hence, by Gauduchon's Plurigenera Theorem \cite{G1} all plurigenera of $(M,I)$ vanish (see Proposition~I.18 and  Proposition~I.19 in \cite{G4} or \cite{V2}). Thus the Kodaira dimension of $(M,I)$ is $-\infty$.

The Kodaira--Enriques classification \cite{BPV} combined with the results in \cite{Kod2,I,Bog} shows that the minimal complex surfaces of Kodaira dimension $-\infty$ are:
$$
\begin{array}{ll}
\left.
\begin{array}{l}
\bullet \quad \mathbb{C}P^2 \\
\bullet \quad \mbox{ruled surfaces}\qquad \qquad
\end{array}
 \right\} &
\mbox{K\"ahler type  ($b_1$ even)} \\
\left.
\begin{array}{l}
\bullet \quad \mbox{Hopf surfaces ($b_2=0$)} \\
\bullet \quad \mbox{Inoue surfaces ($b_2=0$)} \\
\bullet \quad \mbox{surfaces with $b_2>0$}
\end{array}
 \right\} &
\mbox{non-K\"ahler type ($b_1 =1$)}
\end{array}
$$

It is well known that $\mathbb{C}P^2$ is not spin.

Now we show that the Inoue surfaces do not admit metrics with positive scalar curvature by applying Theorem~5.4 in \cite{LM}. More precisely, we use its proof according to which a bundle over a torus with enlargeable fibres is itself enlargeable. There are three types of Inoue surfaces \cite{I}. The surfaces of the first type ($S_M$) are diffeomorphic to 3-torus bundles over a circle \cite{I} and are therefore enlargeable. The second type consists of the surfaces $S^{(+)}_{N,p,q,r,t}$. They are diffeomorphic to bundles over a circle with fibres which are circle bundles over a 2-torus \cite{I}. Thus the fibres are enlargeable and therefore $S^{(+)}_{N,p,q,r,t}$ are enlargeable. Every surface of the third type ($S^{(-)}_{N,p,q,r}$) is double covered by an Inoue surface of second type and is therefore also enlargeable. Hence all Inoue surfaces are enlargeable and by a theorem of Gromov--Lawson \cite{GL,LM} they do not carry metrics with positive scalar curvature. In particular, they cannot have Hermitian metrics with positive conformal scalar curvature.

Finally, let $M$ be a surface of non-K\"ahler type with Kodaira dimension $-\infty$ and $b_2>0$. By Theorem~3 in \cite{Kod1} $b_2^+ = 0$,  $b_2^- = b_2$ and therefore the signature $\sigma(M) <0$. But on a spin manifold with positive scalar curvature the index of the Dirac operator $D$ vanishes and so $\sigma(M) = -8 ind(D) =0$. Thus there are no surfaces with $b_2>0$ which satisfy the assumptions of the proposition.

Hence $(M,I)$ is either a ruled surface or a Hopf surface.  \hfill $\Box $

\vspace{3mm}
This proposition shows that to find the limiting surfaces we need to study the ruled surfaces and the Hopf surfaces.

\section{Ruled surfaces}

\vspace{3mm}
\noindent {\bf Definition \cite{Bea,Har}} A complex surface $M$ is {\it a ruled surface of genus $g$} if it is a holomorphic $\mathbb{C}P^1$-bundle over a compact complex curve $C$ of genus $g$. This is equivalent to $M$ being the projectivization $P(E)$ of some holomorphic vector bundle $E$ of rank 2 over $C$.

\vspace{3mm}
According to a theorem of Yau \cite{Y} every ruled surface carries a K\"ahler metric with positive scalar curvature. Further, Vaisman \cite{V1} has proved that if a complex surface admits a K\"ahler metric, then every locally conformally K\"ahler metric on $M$ is
globally conformally K\"ahler. Thus a locally conformally K\"ahler metric with positive conformal scalar curvature on a ruled surface is globally conformal to a K\"ahler metric with positive scalar curvature. The Lee form of a K\"ahler metric vanishes. Hence Proposition~\ref{prop0} shows that the limiting ruled surfaces are those which admit a holomorphic bundle $S$ such that $S^2 = K$ (i.e., which are spin) and $H^0(M,\mathcal{O}(S^*)) \not = 0$.

Let $M=P(E)$ be a ruled surface and $\pi : P(E) \longrightarrow C$ be the projection. If $L$ is a line bundle over $C$, then obviously $P(E \otimes L) = P(E)$ (the converse is also true: if $P(E') = P(E)$, then $E' = E \otimes L$). Thus, if $e' \in \mathbb{Z}$ has the same parity as $\deg E$, we can represent $M$ as  $P(E')$ with $\deg E' = e'$.

Let $H_E$ denote the tautological line bundle on $P(E)$. Its fibre at $x \in P(E)$ is the line in $(\pi^* E)_x \cong E_{\pi(x)}$ to which $x$ corresponds. Notice that $H_E$ depends on the choice of $E$: $H_{E \otimes L} = H_E \otimes \pi^* L$.

The canonical bundle of $P(E)$ is $K = H_E^2 \otimes \pi^* (\Lambda^2 E^* \otimes K_C)$, where $K_C$ is the canonical bundle of $C$.

We are looking for line bundles $S$ on $P(E)$ such that $S^2 = K$. Since the Picard group of $P(E)$ is generated by $H_E$ and the pull-back of the Picard group of $C$, such an $S$ has the form $S= H_E \otimes\pi^* S_1$, where $S_1$ is a line bundle on $C$ satisfying
\begin{equation}\label{6}
S_1^2 = \Lambda^2 E^* \otimes K_C.
\end{equation}
As $\deg K_C = 2g-2$, such an $S_1$ exists iff $\deg E$ is even. Thus a ruled surface is spin iff  $\deg E$ is even and the number of its spin structures is equal to the number of holomorphic square roots of a line bundle on $C$, i.e., to $2^{2g}$. In fact, the ruled surfaces with even $\deg E$ are homeomorphic to the trivial bundle $C \times S^2$ (and therefore are spin) while those with odd $\deg E$ are homeomorphic to the non-trivial $S^2$-bundle over $C$.

We want also $H^0(M,\mathcal{O}(S^*)) \not = 0$. Since $S^* = H_E^* \otimes \pi^* S_1^*$, by Theorem~5.1 in \cite{BPV} (or Lemma~2.4 in \cite{Har}) $H^0(M,\mathcal{O}(S^*)) = H^0(C,\mathcal{O}(E^* \otimes S_1^*))$. Thus we have to find all rank 2 bundles $E$ and line bundles $S_1$ on $C$ which satisfy \pref{6} and
\begin{equation}\label{7}
H^0(C,\mathcal{O}(E^* \otimes S_1^*)) \not = 0.
\end{equation}

By multiplying by suitable line bundle we can always assume that the bundle $E$ defining the ruled surface has the following property: $H^0(C,\mathcal{O}(E^*)) \not = 0$ but $H^0(C,\mathcal{O}(E^* \otimes L)) = 0$ for each line bundle $L$ with $\deg L <0$.  Such an $E$ is called {\it normalized} \cite{Har}. It may be not unique but $e := \deg E$ does not depend on the particular choice of a normalized $E$ and is therefore an invariant of the ruled surface. (Warning: Our notation $P(E)$ is as in \cite{Bea} and \cite{BPV}. The notation $P(E)$ in \cite{Har} corresponds to our $P(E^*)$.)

From now on we assume that $E$ is normalized. Thus \pref{7} implies that $\deg S_1^* \ge 0$. By \pref{6} we obtain $\deg S_1^* = \frac{1}{2} e + 1 - g$. Hence
\begin{equation}\label{8}
e \ge 2g-2.
\end{equation}

A rank 2 bundle over a complex curve $C$ with $g=0$ (i.e., over $\mathbb{C}P^1$) is decomposable. According to Theorem~2.12 in \cite{Har}, if $g>0$ and a normalized $E$ is indecomposable, then $e \le 2g-2$. Thus \pref{8} shows that in our case an indecomposable $E$ could only occur if $e = 2g-2$. It is proved in Theorem~2.15 in \cite{Har} that if $g=1$, then there exists a unique indecomposable normalized $E$ with $e=0$ and it is the unique non-trivial extension
$$0 \longrightarrow \pmb{1}_C \longrightarrow E \longrightarrow \pmb{1}_C \longrightarrow 0.$$
The same proof can be modified in a straightforward way to show that if $g>0$, then there exists a unique indecomposable normalized $E$ with $e=2g-2$ and it is the unique non-trivial extension
\begin{equation}\label{9}
0 \longrightarrow K_C \longrightarrow E \longrightarrow \pmb{1}_C \longrightarrow 0.
\end{equation}
Thus $\Lambda^2 E \cong K_C$ and therefore by \pref{6} $S_1^2 = \pmb{1}_C$. In particular, $\deg S_1 = 0$. From \pref{9} we obtain an exact sequence
$$0 \longrightarrow S_1^* \longrightarrow E^* \otimes S_1^* \longrightarrow K_C^* \otimes S_1^* \longrightarrow 0.$$
Now \pref{7} implies $H^0(C,\mathcal{O}(S_1^*)) \not = 0$ or $H^0(C,\mathcal{O}(K_C^* \otimes S_1^*)) \not = 0$. If $H^0(C,\mathcal{O}(K_C^* \otimes S_1^*)) \not = 0$, then $0 \le \deg (K_C^* \otimes S_1^*) = 2-2g$. Hence $g=1$, $K_C = \pmb{1}_C$ and $H^0(C,\mathcal{O}(K_C^* \otimes S_1^*)) =  H^0(C,\mathcal{O}(S_1^*))$. This shows that in any case $H^0(C,\mathcal{O}(S_1^*)) \not = 0$ and therefore $S_1 = \pmb{1}_C$ since $\deg S_1 = 0$. Thus $S = H_E \otimes\pi^* \pmb{1}_C = H_E$.

Now let us consider the case of decomposable $E$. Since $E$ is normalized, it has the form $E = L \oplus \pmb{1}_C$, where $\deg L \ge 0$. In particular, $\Lambda^2 E \cong L$ and therefore $e \ge 0$. By \pref{6} this also implies $L = {S_1^*}^2 \otimes K_C$. We have $H^0(C,\mathcal{O}(E^* \otimes S_1^*)) = H^0(C,\mathcal{O}(S_1^*)) \oplus H^0(C,\mathcal{O}(L^* \otimes S_1^*))$. Thus, by \pref{7}, $H^0(C,\mathcal{O}(S_1^*)) \not =0$ or $H^0(C,\mathcal{O}(L^* \otimes S_1^*)) \not = 0$. Let $H^0(C,\mathcal{O}(L^* \otimes S_1^*)) \not = 0$. Then $ 0 \le \deg (L^* \otimes S_1^*) = - \frac{1}{2} e + 1 - g$. Thus $0 \le e \le 2 - 2g$, i.e., $g=0$ or $g=1$. If $g=0$, then $H^0(C,\mathcal{O}(S_1^*)) \not =0$ as $\deg S_1^* \ge 0$. If $g=1$, then $K_C = \pmb{1}_C$ and therefore $L^* \otimes S_1^* = S_1$. So $H^0(C,\mathcal{O}(S_1)) \not =0$, whence $S_1 = \pmb{1}_C$ since $\deg S_1 \le 0$. Thus in any case we have $H^0(C,\mathcal{O}(S_1^*)) \not =0$. This implies, in particular, that either $\deg S_1^* > 0$ (i.e., $e > 2g-2$) or $S_1 = \pmb{1}_C$. In the latter case $e = 2g-2 \ge 0$ and $S = H_E \otimes\pi^* \pmb{1}_C = H_E$.

We summarize the obtained results in the following

\begin{theo}\label{th2}
A compact Hermitian spin surface $(M,h,I,S)$ with even first Betti number is a limiting surface iff $M=P(E)$ is a ruled surface over a complex curve $C$ of genus $g$, $h=c s_0 h_0$, where $c$ is positive constant, $h_0$ is
a  K\"ahler metric with positive scalar curvature $s_0$, and $E$ and $S$ are given by:

1) $g >0$, $E$ is the unique non-trivial extension of $\pmb{1}_C$ by $K_C$, $S = H_E$. In this case $E$ is indecomposable.

2) $g >0$, $E = K_C \oplus \pmb{1}_C$ , $S = H_E$.

3) $g$ is arbitrary, $E= N^2 \otimes K_C \oplus \pmb{1}_C$, $S=H_E \otimes \pi^* N^*$, where $N$ is a line bundle over $C$ with $\deg N >0$ and $H^0(C,\mathcal{O}(N)) \not = 0$.
\end{theo}

We note that all bundles $E$ in this theorem are normalized. In fact, for all ruled surfaces appearing in it there exists a unique normalized $E$. There is also no repetition of pairs $(E,S)$ and thus a ruled surface appears more than once if it has more than one spin structure with which it is limiting. The surfaces in cases 1) and 2) are limiting  with only one of their spin structures. At the other extreme, the surfaces $P(N^2 \otimes K_C \oplus \pmb{1}_C)$ with $\deg N > g$ are limiting with all of their $2^{2g}$ spin structures. Furthermore, every ruled surface of genus $0$ which is spin appears in the above list, i.e., the even Hirzebruch surfaces with their unique spin structure are limiting.

Finally, the results in \cite{F2} show that the K\"ahler limiting surfaces, i.e., the limiting manifolds for Kirchberg's estimate \pref{2}, are $\mathbb{C}P^1 \times \mathbb{C}P^1$ (which comes from 3) with $g=0$ and $N = K_{\mathbb{C}P^1}^{-1/2}$) and $C \times \mathbb{C}P^1$, where $C$ is a torus (this comes from 2) with $g=1$). Thus there exist many more limiting manifolds for Theorem~\ref{th1} than in the K\"ahler case. This is not surprising since the existence of a K\"ahler metric with positive constant scalar curvature implies strong restrictions on the automorphism group of the complex surface.

\section{Hopf surfaces}

\noindent {\bf Definition \cite{Kod2}} {\it A Hopf surface} is a compact complex surface whose universal cover is $\mathbb{C}^2 \backslash 0$. {\it A primary Hopf surface} is a Hopf surface whose fundamental group is $\mathbb{Z}$.

\vspace{3mm}
It was shown in \cite{Kod2} that  every Hopf surface is finitely covered by a primary one and the primary Hopf surfaces are of the form $M = (\mathbb{C}^2 \backslash 0)/<\gamma>$, where
$$\gamma(z_1,z_2)=(\alpha_1 z_1 + \lambda z_2^m,\alpha_2 z_2),$$
$$\alpha_1, \alpha_2 \in \mathbb{C}, \quad 0< |\alpha_1| \le |\alpha_2| <1, \quad m \in \mathbb{N}, \quad \lambda \in \mathbb{C}, \quad \lambda (\alpha_1 - \alpha_2^m) =0.$$
A primary Hopf surface is said to be of class 1 if $\lambda =0$ (and we denote it by $M_{\alpha_1, \alpha_2}$) and of class 0 if $\lambda \not =0$, $\alpha_1 = \alpha_2^m$ (we denote it by $M_{\alpha_1, \alpha_2,m,\lambda}$). It is clear that $M_{\alpha_1, \alpha_2, m, \lambda}$ and $M_{\alpha_1, \alpha_2, m, \mu}$ are biholomorphically equivalent for arbitrary $\lambda \not =0$, $\mu \not =0$. Furthermore, a primary Hopf surface is diffeomorphic to $S^1 \times S^3$ \cite{Kod2}. An explicit diffeomorphism can be obtained in the following way:

Let $c = 1 + \left| \frac{\lambda}{\alpha_1 \ln |\alpha_2|} \right|^{\frac{1}{m}}$ and  $t:\mathbb{C}^2 \backslash 0 \longrightarrow \mathbb{R}$ be the function defined by the equation
$$|\alpha_1|^{-2t(z)} |z_1 - \frac{\lambda}{\alpha_1} t(z) z_2^m|^2 + c^2 |\alpha_2|^{-2t(z)} |z_2|^2 = 1.$$
Define $\widetilde{g} : \mathbb{C}^2 \backslash 0 \longrightarrow S^1 \times S^3 \subset \mathbb{C} \times \mathbb{C}^2$ through
$$z \mapsto (e^{2\pi i t(z)},\alpha_1^{-t(z)} (z_1 - \frac{\lambda}{\alpha_1} t(z) z_2^m), c \alpha_2^{-t(z)} z_2).$$
Here we have fixed $\varphi_1, \varphi_2 \in \mathbb{R}$ so that $\alpha_j = |\alpha_j| e^{i\varphi_j}$ and $\alpha_j^r := |\alpha_j|^r e^{ir\varphi_j}$ for $r \in \mathbb{R}$, $j=1,2$. Then $\widetilde{g}$ is $\gamma$-invariant and covers a diffeomorphism $g:M \longrightarrow S^1 \times S^3$, i.e., $\widetilde{g} = g \circ \pi$, where $\pi : \mathbb{C}^2 \backslash 0 \longrightarrow M$ is the projection. Notice also that $t(\gamma(z)) = t(z) +1$.

Let $\omega_0 = \frac{1}{2\pi} (-ydx + xdy)$ be the form generating $H^1(S^1,\mathbb{Z})$. Then $\omega = g^* \omega_0$ generates $H^1(M,\mathbb{Z})$. We have $\pi^* \omega =  \widetilde{g}^* \omega_0 = dt$.

It is well known \cite{Kod3} that on a Hopf surface $M$ $H^1(M,\mathcal{O}^*) = H^1(M,\mathbb{C}^*)$ and if $M$ is furthermore primary, then $H^1(M,\mathbb{C}^*) \cong \mathbb{C}^*$. This isomorphism can be realized by an explicit construction:

\begin{lem}\label{lem2}
Let $M$ be a primary Hopf surface. Then there exists an isomorphism  $\mathbb{C}^* \cong H^1(M,{\mathcal O}^*)$ assigning to $s \in \mathbb{C}^*$ a bundle $L(s) \in H^1(M,\mathbb{C}^*) = H^1(M,{\mathcal O}^*)$ with the following properties:

1) $K = L(\frac{1}{\alpha_1 \alpha_2})$

2) $E([\mu \omega]) = L(e^{2\pi i \mu})$, $\mu \in \mathbb{C}$

3) $H^0(M,L(s)) \not = 0$ iff $s = \alpha_1^{n_1} \alpha_2^{n_2}$, $n_1,n_2 \in \mathbb{N}_0$.
\end{lem}

\noindent {\it Proof:}
The proof is a straightforward generalization of those given in \cite{KS,Bes4} for $M_{\alpha,\alpha}$:

The bundle $L(s)$ is defined by
$$L(s) = ((\mathbb{C}^2 \backslash 0) \times \mathbb{C})/<\gamma_s>,$$
where $\gamma_s (z,\xi) = (\gamma (z),s\xi)$.
Conversely, given $L \in H^1(M,\mathbb{C}^*)$, then $L \cong L(s)$, where $s \in \mathbb{C}^*$ is obtained in the following way:

Since $H^1(\mathbb{C}^2 \backslash 0,\mathbb{C}^*)$ is trivial, we can fix a trivialization $\pi^* L \cong (\mathbb{C}^2 \backslash 0) \times \mathbb{C}$. Let $0 \not = l \in L_{\pi(z)}$. Then $(\pi^* l)_z = (z, \xi_0)$, $(\pi^* l)_{\gamma(z)} = (\gamma(z), \xi_1)$ and $s = \frac{\xi_1}{\xi_0}$.

Let us apply this construction to $K$. We have $\pi^* K = K_{\mathbb{C}^2 \backslash 0}$ which is trivialized by $dz_1 \wedge dz_2$. Let $\varphi \in K_{\pi(z)}$, $(\pi^* \varphi)_z = \xi_0 (dz_1 \wedge dz_2)_z$, $(\pi^* \varphi)_{\gamma(z)} = \xi_1 (dz_1 \wedge dz_2)_{\gamma(z)}$. We have $\gamma^* ((\pi^* \varphi)_{\gamma(z)}) = (\pi^* \varphi)_z$ which yields $\xi_1 \alpha_1 \alpha_2 (dz_1 \wedge dz_2)_z = \xi_0 (dz_1 \wedge dz_2)_z$. Hence $\frac{\xi_1}{\xi_0} = \frac{1}{\alpha_1 \alpha_2}$, i.e., $K = L(\frac{1}{\alpha_1 \alpha_2})$.

Now consider a closed 1-form $\theta$ on $M$. Then $\pi^* \theta = df$ for some function $f$ on $\mathbb{C}^2 \backslash 0$. We have $\pi^* E([\theta]) = E([\pi^* \theta])$ and using this it is straightforward to see that $E([\theta]) = L(e^{2\pi i (f(\gamma(z)) - f(z))})$ for arbitrary $z \in \mathbb{C}^2 \backslash 0$. Since $\pi^* \omega = dt$, we obtain $E([\mu \omega]) = L(e^{2\pi i \mu (t(\gamma(z)) - t(z))}) = L(e^{2\pi i \mu})$.

A holomorphic section of $L(s)$ is pulled back to a holomorphic section of $(\mathbb{C}^2 \backslash 0) \times \mathbb{C}$ which is identified with a holomorphic function $\varphi : \mathbb{C}^2 \backslash 0 \longrightarrow \mathbb{C}$ such that
\begin{equation}\label{41}
\varphi (\gamma(z)) = s \varphi (z).
\end{equation}
By Hartogs theorem $\varphi $ extends to a holomorphic function on $\mathbb{C}^2$. It can be  written as a power series $\varphi (z) = \sum_{k_1,k_2 = 0}^{\infty} C_{k_1k_2} z_1^{k_1} z_2^{k_2}$ and \pref{41} gives equations for $C_{k_1k_2}$ which have non-zero solutions only if $s = \alpha_1^{n_1} \alpha_2^{n_2}$ with $n_1,n_2 \in \mathbb{N}_0$. \hfill $\Box$

\vspace{3mm}
\noindent {\bf Remark:} For an arbitrary Hopf surface one can define in a similar way an explicit isomorphism $H^1(M,\mathbb{C}^*) \cong Hom(\pi_1 (M), \mathbb{C}^*)$.

\vspace{3mm}
This lemma shows in particular that on a primary Hopf surface there are only two spin structures $S$. The corresponding $S^*$ are $L(\pm \sqrt{\alpha_1 \alpha_2})$, where $\sqrt{\alpha_1 \alpha_2}$ is a fixed square root of $\alpha_1 \alpha_2$. Thus Proposition~\ref{prop0} and Lemma~\ref{lem2} yield

\begin{theo}\label{th3}
Let $(M,h,I)$ be a primary Hopf surface with locally conformally K\"ahler metric $h$ with positive conformal scalar curvature. Let $\mu \in \mathbb{R}^*$ be determined by $[\theta]= [\mu \omega]$, where $\theta$ is the Lee form of $h$. Then there exists a spin structure $S$ with which $(M,h,I,S)$ is a limiting surface iff there exist $n_1,n_2 \in \mathbb{N}_0$ such that $\alpha_1^{2n_1-1} \alpha_2^{2n_2-1} =e^{-2\mu}$.
\end{theo}

Recall that there is a particular class of locally conformally K\"ahler manifolds, the generalized Hopf (or Vaisman) manifolds. These are Hermitian manifolds whose Lee form is parallel with respect to the Levi-Civita connection. It was proved in \cite{Bel} that a Hopf surface admits a Vaisman metric only if it is of class 1. For primary Hopf surfaces of class 1 an explicit Vaisman metric was constructed by Gauduchon--Ornea \cite{GO}. It is defined by its K\"ahler form $\Omega$ through
$$\pi^* \Omega = - \frac{e^{-\ln |\alpha_1 \alpha_2| \, t}}{4} dd^c e^{\ln |\alpha_1 \alpha_2| \, t}.$$
Since $- \frac{1}{4} dd^c e^{\ln |\alpha_1 \alpha_2| \, t}$ is the K\"ahler form of a K\"ahler metric on $\mathbb{C}^2 \backslash 0$, for the Lee form we obtain
$$\pi^* \theta = -\ln |\alpha_1 \alpha_2| \, dt = -\ln |\alpha_1 \alpha_2| \, \pi^* \omega $$
and therefore $\theta = -\ln |\alpha_1 \alpha_2| \, \omega $.

According to a result of Tsukada \cite{Ts}, given a compact complex manifold $(M,I)$ with Vaisman metric with Lee form $\theta$, then there exists a locally conformally K\"ahler metric with Lee form $\theta_1$ such that $[\theta_1] = [\mu \theta]$ iff $\mu >0$. This implies, in particular, that for primary Hopf surfaces of class 1 we have $\mu >0$ in Theorem~\ref{th3}.

Let us now recall the construction which gives the existence part of Tsukada's result.

Let $(M,h,I)$ be a Hermitian surface with Vaisman metric $h$ with Lee form $\theta$ and conformal scalar curvature $k$. Let $\nu \in \mathbb{R}$. Define $h_\nu = h + \nu (\theta \otimes \theta + I\theta \otimes I\theta)$.  Then it is straightforward to see that $h_\nu$ is a Vaisman metric iff $\nu > - \frac{1}{|\theta|^2}$ (otherwise $h_\nu$ is not positive) and its Lee form and conformal scalar curvature are $\theta_\nu = (1+ \nu |\theta|^2)\theta$, $k_\nu = k - 2\nu |\theta|^4$.

We apply this construction to the metrics of Gauduchon--Ornea. They have
$$|\theta| = 2, \quad \pi^* k = - \frac{24\ln |\frac{\alpha_1}{\alpha_2}| \, (\ln |\alpha_1| \, |z_1|^2 - \ln |\alpha_2| \, |z_2|^2)}{\ln |\alpha_1 \alpha_2| \, (\ln |\alpha_1| \, |z_1|^2 + \ln |\alpha_2| \, |z_2|^2)}.$$
Now it is straightforward to see that $k_\nu > 0$ iff
$$\nu < - \frac{3\ln |\frac{\alpha_1}{\alpha_2}|}{4\ln |\alpha_1 \alpha_2|}.$$
This and $\nu > - \frac{1}{4}$ imply $2\ln |\alpha_2| < \ln |\alpha_1|$, i.e., $|\alpha_2|^2 < |\alpha_1| \le |\alpha_2|$. We have
$$\theta_\nu = (1+4\nu)\theta = - (1+4\nu)\ln |\alpha_1 \alpha_2| \, \omega ,$$
so the corresponding $\mu$ is
$$\mu = - (1+4\nu)\ln |\alpha_1 \alpha_2| < 2(\ln |\alpha_1| - 2\ln |\alpha_2|).$$
Therefore
$$1> \alpha_1^{2n_1-1} \alpha_2^{2n_2-1} =e^{-2\mu} > \frac{|\alpha_2|^8}{|\alpha_1|^4},$$
which yields the following possibilities for $n_1$ and $n_2$:
$$n_1 =0, \quad n_2 =2; \qquad n_1 =1, \quad n_2 =0; \qquad n_1 =1, \quad n_2 =1; \qquad n_1 =2, \quad n_2 =0.$$
Thus for the following values of $\alpha_1$ and $\alpha_2$ the corresponding Hopf surfaces of class 1 admit metrics with which they become limiting surfaces:
$$
\begin{array}{lll}
\alpha_1 = |\alpha_1| e^{i\varphi}, & \alpha_2 = |\alpha_2| e^{\frac{i(\varphi + 2j\pi)}{3}}, \, j=0,1,2, & |\alpha_2|^{\frac{5}{3}} < |\alpha_1| \le |\alpha_2|; \\
\alpha_1 = |\alpha_1| e^{i\varphi}, & \alpha_2 = |\alpha_2| e^{i\varphi}, & |\alpha_2|^{\frac{9}{5}} < |\alpha_1| <|\alpha_2|; \\
\alpha_1 = |\alpha_1| e^{i\varphi}, & \alpha_2 = |\alpha_2| e^{-i\varphi}, & |\alpha_2|^{\frac{7}{5}} < |\alpha_1| \le |\alpha_2|;  \\
\alpha_1 = |\alpha_1| e^{\frac{i(\varphi + 2j\pi)}{3}}, \, j=0,1,2, & \alpha_2 = |\alpha_2| e^{i\varphi}, & |\alpha_2|^{\frac{9}{7}} < |\alpha_1| \le |\alpha_2|.
\end{array}
$$

On the other hand, there are infinitely many Hopf surfaces which do not admit such metrics, for example those for which $\alpha_1^{2n_1-1} \alpha_2^{2n_2-1}$ is never real. And there are also pairs $\alpha_1$, $\alpha_2$ such that there exist $\mu$, $n_1$, $n_2$ with  $\alpha_1^{2n_1-1} \alpha_2^{2n_2-1} = e^{-2\mu}$ but there is no known locally conformally K\"ahler metric with $k>0$ and $[\theta]= [\mu \omega]$ on them.

Thus it is clear that to complete the classification of the limiting surfaces one needs to answer the following question:

Given a Hopf surface, which are the values of $\mu \in \mathbb{R}$ so that there exists a locally conformally K\"ahler metric with positive conformal scalar curvature and Lee form $\mu \omega$?

The only restrictions we know are that $\mu \not = 0$ (since the Hopf surfaces do not admit K\"ahler metrics) and that $\mu >0$ in case of a Hopf surface of class 1.

\vspace{10mm}
\noindent
Bogdan Alexandrov \\
Humboldt Univesit\"at zu Berlin \\
Institut f\"ur Mathemathik \\
Sitz: Rudower Chaussee 25  \\
10099 Berlin \\
{\tt e-mail: \quad bogdan@mathematik.hu-berlin.de}

\end{document}